\documentclass[11pt, a4paper]{amsart}
\usepackage{color}
\def\RR{{\mathbb {R}}}

\def\ep{\varepsilon}

\def\di{\displaystyle}

\def\R{{\mathbb {R}}}
\def\uep{u_{\ep}}
\def\wep{w_{\ep}}

\def\K{{\mathcal{K}}}
\def\cald{{\mathcal{D}}}
\def\calr{{\mathcal{R}}}
\def\J{{\mathcal{J}}}
\def\H{{\mathcal{H}}^{N-1}}

\def\pint{\operatorname {--\!\!\!\!\!\int\!\!\!\!\!--}}
\def\pep{P_\ep}
\def\jep{{\mathcal J}_\ep}

\def\wp{W^{1,p}(\Omega)}

\newtheorem{teo}{Theorem}[section]
\newtheorem{lema}{Lemma}[section]
\newtheorem{defi}{Definition}[section]

\newtheorem{corol}{Corollary}[section]
\theoremstyle{definition}

\newtheorem*{ack}{Acknowledgements}
\renewcommand{\theequation}{\arabic{section}.\arabic{equation}}

\begin{document}

\title[An optimization problem with volume constrain]
{An optimization problem with volume constrain for a degenerate
quasilinear operator}

\author[J. Fern\'andez Bonder, S. Mart\'{i}nez \& N. Wolanski]
{Juli\'an Fern\'andez Bonder, Sandra Mart\'{\i}nez and Noemi Wolanski}

\address{Departamento  de Matem\'atica, FCEyN
\hfill\break\indent UBA (1428) Buenos Aires, Argentina.}
\email{{\tt jfbonder@dm.uba.ar \\ smartin@dm.uba.ar \\ wolanski@dm.uba.ar}
\hfill\break\indent {\em Web-page:}
{\tt http://mate.dm.uba.ar/$\sim$jfbonder}, 
\hfill\break\indent\indent\indent\indent\indent\ {\tt http://mate.dm.uba.ar/$\sim$wolanski}}

\thanks{Supported by ANPCyT PICT No.
03-13719 and 03-10608, UBA X052 and X066 and
Fundaci\'on Antorchas 13900-5. J. Fern\'andez Bonder and N. Wolanski are members
of CONICET}


\begin{abstract}
We consider the  optimization problem of minimizing
$\int_{\Omega}|\nabla u|^p\, dx$ with a constrain on the volume of
$\{u>0\}$. We consider a penalization problem, and we prove that
for small values of the penalization parameter, the constrained
volume is attained. In this way we prove that every solution $u$
is locally Lipschitz continuous and that the free boundary,
$\partial\{u>0\}\cap \Omega$, is smooth.
\end{abstract}

\maketitle

\section{Introduction}

In the seminal paper \cite{AAC}, Aguilera, Alt and Caffarelli
study an optimal design problem with a volume constrain by
introducing a penalization term in the energy functional (the
Dirichlet integral) and minimizing without the volume constrain.
For fixed values of the penalization parameter, the penalized
functional is very similar to the one considered in the paper
\cite{AC}. So that, regularity results for minimizers of the
penalized problem follow almost without change as in \cite{AC}.
The main result in \cite{AAC} that makes this method so useful is
that the right volume is already attained for small values of the
penalization parameter. In this way, all the regularity results
apply to the solution of the optimal design problem.

This method has been applied to other problems with similar
success. In all those cases, the differential equation satisfied
by the minimizers is nondegenerate, uniformly elliptic. See, for
instance, \cite{ACS, FBRW, L, Te}.

In this article we want to show that the same kind of results can
be obtained for some nonlinear degenerate or singular elliptic
equations. As an example, we study here the following problem
which is a generalization of the one in \cite{AAC} for $1<
p<\infty$:

We take $\Omega$ a smooth bounded domain in $\R^N$ and $\varphi_0
\in W^{1,p}(\Omega)$, a Dirichlet datum, with $\varphi_0\geq
c_0>0$ in $\bar{A}$, where A is a nonempty relatively open subset
of $\partial\Omega$ such that $A\cap\partial\Omega$ is $C^2$. Let
$$
\K_\alpha = \{u\in W^{1,p}(\Omega)\, / \, |\{u>0\}| = \alpha, \,
u=\varphi_0\ \mbox{ on }\partial\Omega\}.
$$

Our problem is to minimize $\J(u)=\int_\Omega|\nabla u|^p\,dx$ in
$\K_\alpha$.

Problems similar to the one considered here appear in shape
optimization. For instance, in optimization of torsional rigidity
\cite{L}, insulation of pipelines for hot liquids \cite{Fl},
minimization of the current leakage from insulated wires and
coaxial cables \cite{A}, minimization of the capacity of
condensers and resistors, etc.

Although the existence of a minimizer is not difficult to
establish by variational techniques, the regularity properties of
such minimizers and their free boundaries $\partial\{u>0\}$, are
not easy to obtain since it is hard to make enough volume
preserving perturbations without the previous knowledge of the
regularity of $\partial\{u>0\}$.

In order to solve our original problem in a way that allows us to
perform non volume preserving perturbations we consider instead
the following penalized problem: We let
$$
\K=\{u\in W^{1,p}(\Omega)\, / \, u = \varphi_0\ \mbox{ on }
\partial\Omega\}
$$
and
\begin{equation}\label{jep}
\jep(u) = \int_\Omega|\nabla u|^p\, dx + F_\ep(|\{u>0\}|),
\end{equation}
where
$$
F_\ep(s) = \begin{cases}
\ep(s-\alpha) & \mbox{ if } s<\alpha\\
\frac1\ep(s-\alpha) & \mbox{ if } s\ge\alpha.
\end{cases}
$$

Then, the penalized problem is
\begin{equation*}\label{Pep}
\tag{$P_\ep$} \mbox{Find } \uep\in \K \quad \mbox{such that} \quad
\J_\ep(\uep) = \inf_{v\in \K} \J_\ep(v).
\end{equation*}

The existence of minimizers follows easily by direct minimization.
Their regularity and the regularity of their free boundaries
$\partial \{ \uep >0 \}$ follow as in \cite{DP} where a very
similar problem was studied, namely, to minimize
\begin{equation}\label{dp}
\overline{\J}_\lambda(v) = \int_{\Omega} |\nabla v|^p \, dx +
\lambda^p |\{v>0\}| ,
\end{equation}
where $\lambda>0$ is a constant. In particular, $\uep$ is a solution of the
following free boundary problem
$$
\begin{cases}
\Delta_p u  = 0 &\mbox{ in } \{u>0\}\cap \Omega,\\
\di\frac{\partial u}{\partial\nu} = \lambda_\ep & \mbox{ on }
\partial\{u>0\}\cap\Omega,
\end{cases}
$$
where $\lambda_\ep$ is a positive constant and $\Delta_p u = {\rm
div}(|\nabla u|^{p-2}\nabla u)$ is the $p-$laplacian.

In \cite{FBRW} the authors study a problem closely related to
\cite{AAC}. The problem in \cite{FBRW} is to minimize the best
Sobolev trace constant from $H^1(\Omega)$ into
$L^q(\partial\Omega)$ for subcritical $q$, among functions that
vanish in a set of fixed measure. We will sometimes refer to some
of the proofs in \cite{FBRW} for the different treatment of the
penalization term (which is piecewise linear in the measure of the
positivity set) with respect to \cite{AC} and \cite{DP} where the
function is linear in the measure.

As in \cite{AAC}, the reason why this penalization method is so
useful is that there is no need to pass to the limit in the
penalization parameter $\ep$ for which uniform, in $\ep$,
regularity estimates would be needed. In fact, we show that for
small values of $\ep$ the right volume is already attained. This
is, $|\{\uep>0\}|=\alpha$ for $\ep$ small. It is at this point
where the main changes have to be made since the perturbations
used in \cite{AAC} and \cite{FBRW} make strong use of the
linearity of the underlying equation.

In particular, the fact that, for small $\ep$, any minimizer of
$\J_\ep$ satisfies $|\{\uep>0\}|=\alpha$ implies that any
minimizer of our original optimization problem is also a minimizer
of $\J_\ep$ so that it is locally Lipschitz continuous with smooth
free boundary.

\bigskip

We include at the end of the paper a couple of appendices where
some properties of $p-$sub\-harmonic functions are established. We
use these results in Section 2. We believe that these results have
independent interest.

The paper is organized as follows: In Section 2 we begin our
analysis of problem \eqref{Pep} for fixed $\ep$. First we prove
the existence of a minimizer, local Lipschitz regularity and
nondegeneracy near the free boundary (Theorem \ref{existencia})
and with these results we have the regularity of the free boundary
by adapting the results of \cite{DP}.

The  main results of this paper appear in Section 3 where we prove
that for small values of $\ep$ we recover our original
optimization problem.

The appendices are included at the end of the paper.


\section{The penalized problem}

In this section we look for minimizers of the  functional $\J_\ep$
and a representation theorem for solutions of $\J_\ep$ as in
\cite{AC} Theorem 4.5.

Observe that a solution to $(\pep)$ satisfies that
$$
\Delta_p u = 0 \quad \mbox{in } \{u>0\}^\circ.
$$
In fact, let $B$ be a ball such that $u>0$ in $B$. Let $v$ be the
solution to
$$
\Delta_p v=0\quad\mbox{in }B, \qquad v=u \quad \mbox{on } \partial
B.
$$
Let $\bar v\in W^{1,p}(\Omega)$, $\bar v(x)=v(x)$ for $x\in B$,
$\bar v(x)=u(x)$ if $x\notin B$. Then, $\bar v\in \K$ so that
\begin{equation}\label{primeracota}
0 \le \int_\Omega |\nabla\bar v|^p\, dx-\int_\Omega |\nabla u|^p\,
dx + F_\ep(|\{\bar v>0\}|) - F_\ep(|\{u>0\}|)= \int_B |\nabla v|^p
-|\nabla u|^p\, dx,
\end{equation}
and (see \cite{DP}, Section 3),
\begin{align}\label{desigualdadDP1}
&\int_B |\nabla v|^p-|\nabla u|^p\, dx \le -c\int_B
|\nabla(v-u)|^p\, dx \qquad \mbox{if } p\ge2,\\
\label{desigualdadDP2} &\int_B |\nabla v|^p - |\nabla u|^p\, dx
\le -c \int_B |\nabla (v-u)|^2 \big(|\nabla v| + |\nabla u|
\big)^{p-2}\, dx \qquad \mbox{if }1<p\le2,
\end{align}
where $c$ is a positive constant that depends on $p$. In any case,
combining \eqref{primeracota} and \eqref{desigualdadDP1} --
\eqref{desigualdadDP2} we get $|\nabla(v-u)|=0$ in $B$. Thus, $u =
v$ in $B$. So that, $\Delta_p u=0$ in $B$.

\bigskip

We begin by discussing the existence of extremals.

\begin{teo}\label{existencia}
Let $\Omega\subset\R^N$ be bounded and $1<p<\infty$. Then there
exists a solution to the problem \eqref{Pep}. Moreover, any such
solution $\uep$ has the following properties:
\begin{enumerate}
\item $\uep$ is locally Lipschitz continuous in $\Omega$.

\item For every $D\subset\subset\Omega$, there exist constants $C,
c>0$ such that for every $x\in D\cap \{\uep>0\}$,
$$
c\, {\rm dist}(x,\partial\{\uep>0\})\le \uep(x)\le C\, {\rm
dist}(x,\partial\{\uep>0\}).
$$

\item For every $D\subset\subset\Omega$, there exists a constant
$c>0$ such that for $x\in \partial\{u>0\}$ and $B_r(x)\subset D$,
$$
c\le \frac{|B_r(x)\cap\{\uep>0\}|}{|B_r(x)|}\le 1-c.
$$
\end{enumerate}
The constants may depend on $\ep$.
\end{teo}

\begin{proof}
The proof of existence is standard. We state it here for the
reader's convenience.

Take $u_0$ with $|\{u_0>0\}|\leq \alpha$, then $\J_{\ep}(u_0)\leq
C$ (uniformly in $\ep$), also $\J_{\ep}\geq-\alpha$. Therefore a
minimizing sequence $(u_n)\subset\K$ exists. Then $\J_\ep(u_n)$ is
bounded, so $\|\nabla u_n\|_{p}\le C$. As  $u_n = \varphi_0$ in
$\partial\Omega$, there exists a subsequence (that we still call
$u_n$) and a function $\uep\in W^{1,p}(\Omega)$ such that
\begin{align*}
& u_n \rightharpoonup\uep \quad \mbox{weakly in } W^{1,p}(\Omega),\\
& u_n \to \uep \quad \mbox{a.e. } \Omega.
\end{align*}
Thus,
\begin{align*}
& \uep=\varphi_0\quad\mbox{on}\quad \partial\Omega,\\
& |\{\uep>0\}| \le \liminf_{n\to\infty}|\{u_n>0\}| \quad\mbox{ and}\\
& \int_{\Omega}|\nabla \uep|^p\, dx \le \liminf_{n\to\infty}
\int_{\Omega}|\nabla u_n|^p\, dx.
\end{align*}
Hence $\uep\in \K$ and
$$
\J_\ep(\uep)\le \liminf_{n\to\infty}\J_\ep(u_n) = \inf_{v\in\K} \J_\ep(v),
$$
therefore $\uep$ is a minimizer of $\J_\ep$ in $\K$.

The proof of (1), (2) and (3) follow as Theorem 3.3, Lemma 4.2 and
Theorem 4.4 in \cite{DP}. The only difference being that the
functional they analyze is linear in $|\{\uep>0\}|$ and ours is
piecewise linear. The different treatment of this term is similar
to the one in \cite{FBRW}.
\end{proof}

From now on we denote by $u$ instead of $\uep$  a solution to
\eqref{Pep}.

\begin{lema}\label{medida1} Let $u\in \K$ be a solution to \eqref{Pep}.
Then $u$ satisfies for every $\varphi\in C_0^\infty(\Omega)$ such
that ${\rm supp}(\varphi)\subset\{u>0\}$,
\begin{equation}\label{ecuacion}
\int_\Omega|\nabla u|^{p-2}\nabla u\nabla \varphi=0.
\end{equation}
Moreover, the application
$$
\lambda(\varphi):= -\int_\Omega|\nabla u|^{p-2}\nabla u\nabla \varphi\,dx
$$
from $C_0^\infty(\Omega)$ into $\R$ defines a nonnegative Radon
measure  with support on $\Omega\cap\partial\{u>0\}$.
\end{lema}

\begin{proof} See Theorem 5.1 in \cite{DP}
\end{proof}

\begin{teo}[Representation Theorem] \label{DP-teo5.2} Let $u\in \K$ be a
solution to \eqref{Pep}. Then,
\begin{enumerate}
\item $\H( D\cap\partial\{u>0\})<\infty$ for every
$D\subset\subset\Omega$.

\item There exists a Borel function $q_u$ such that
$$
\Delta_p u= q_u \,\H \lfloor \partial\{u>0\}.
$$

\item For $D\subset\subset\Omega$ there are constants $0<c\le
C<\infty$ depending on $N, \Omega, D$ and $\ep$ such that for
$B_r(x)\subset D$ and $x\in \partial\{u>0\}$,
$$
c\le q_u(x)\le C,\quad c\,r^{N-1}\le
\H(B_r(x)\cap\partial\{u>0\})\le C\,r^{N-1}.
$$

\item for $\H$--a.e. $x_0\in \partial_{\rm red}\{u>0\}$,
$$
u(x_0+x) = q_u(x_0)  (x\cdot \nu(x_0))^- + o(|x|) \quad
\mbox{for}\quad x\to 0
$$
with $\nu(x_0)$ the outward unit normal de $\partial\{u>0\}$ in
the measure theoretic sense.

\item $\H(\partial\{u>0\}\setminus\partial_{\rm red}\{u>0\})=0.$
\end{enumerate}
\end{teo}

\begin{proof}
The proof of (1), (2) and (3) follow exactly as that of Theorem
4.5 in \cite{AC}.

Observe that $D\cap\partial\{u>0\}$ has finite perimeter, thus,
the reduce boundary $\partial_{\rm red}\{u>0\}$ is defined as well
as the measure theoretic normal $\nu(x)$ for $x\in\partial_{\rm
red}\{u>0\}$ (see \cite{F}). For the proof of (4) see Theorem 5.5
in \cite{DP}.

Finally, (5) is a consequence of Theorem \ref{existencia} and (3)
(see \cite{F}).
\end{proof}

\begin{teo}\label{const}
Let $u\in \K$ be a solution to \eqref{Pep} and $q_{u}$ the
function in Theorem \ref{DP-teo5.2}. Then there exists a constant
$\lambda_u$ such that
\begin{align}
& \limsup_{\stackrel{x\to x_0}{u(x)>0}} |\nabla u(x)| =
\lambda_u,\qquad \mbox{for every } x_0\in
\Omega\cap\partial\{u>0\}\label{limsup}\\
& q_{u}(x_0)=\lambda_u,\qquad \H -\mbox{a.e } x_0\in
\Omega\cap\partial\{u>0\}. \label{qu}
\end{align}
Moreover, if $B$ is a ball contained in $\{u=0\}$ touching the
boundary $\partial\{u>0\}$ at $x_0$. Then
\begin{equation}\label{limsup.dist}
\hskip-5cm \limsup_{\stackrel{x\to x_0}{u(x)>0}} \frac{u(x)}{{\rm
dist}(x,B)}=\lambda_u.
\end{equation}
\end{teo}

To prove this theorem, we  have to prove first the following lemma,

\begin{lema}\label{blow}
Let $x_0, x_1\in \partial\{u>0\}$ and $\rho_k\to 0^+$. For $i=0,1$
let $x_{i,k}\to x_i$ with $u(x_{i,k})=0$ such that
$B_{\rho_k}(x_{i,k}) \subset \Omega$ and such that the blow-up
sequence
$$
u_{i,k}(x) = \frac{1}{\rho_k} u(x_{i,k} + \rho_k x)
$$
has a limit $u_i(x)=\lambda_i ( x\cdot \nu_i)^-$, with
$0<\lambda_i<\infty$ and $\nu_i$ a unit vector. Then $\lambda_0 =
\lambda_1$.
\end{lema}
\begin{proof}
Assume that $\lambda_1<\lambda_0$ then we will perturb the
minimizer $u$ near $x_0$ and $x_1$ and get an admissible function
with less energy. To this end, we take a nonnegative $C_0^\infty$
function $\phi$ supported in the unit interval. For $k$ large,
define
$$
\tau_k(x) = \begin{cases}
x + \rho_k^2\phi\Big(\frac{\di |x-x_{0,k}|}{\di\rho_k}\Big)\nu_0 &
\mbox{for } x\in B_{\rho_k}(x_{0,k}),\\
\\
x - \rho_k^2\phi\Big(\frac{\di |x-x_{1,k}|}{\di\rho_k}\Big)\nu_1 &
\mbox{for } x\in B_{\rho_k}(x_{1,k}),\\
\\
x & \mbox{elsewhere},
\end{cases}
$$
which is a diffeomorphism if $k$ is big enough. Now let
$$
v_k(x) = u(\tau_k^{-1}(x)),
$$
that are admissible functions. Let us also define
\begin{equation}\label{eta}
\eta_i(y) = (-1)^i\phi(|y|)\nu_i.
\end{equation}
We have
\begin{equation}\label{estimacion.medida}
F_\ep(|\{v_k>0\}|)-F_\ep(|\{u>0\}|) = o(\rho_k^{N+1}).
\end{equation}

To estimate the other term in $\jep$ we make a change of variables
and then
\begin{align*}
\rho_k^{-N}\int_{B_{\rho_k}(x_i)}&(|\nabla v_k|^p-|\nabla
u_{\ep}|^p)\ dx \\
&= \int_{B_{1}(0)\cap \{u_{i,k}>0\}}\rho_k \Big[|\nabla u_{i,k}|^p
\,\mbox{div}(\eta_i)- p \,|\nabla u_{i,k}|^{p-2} (\nabla
u_{i,k})^t D\eta_i \nabla u_{i,k}\Big] + o(\rho_k)\ dy.
\end{align*}

On the other hand, by Lemma \ref{propblowup}, we have
\begin{align*}
& B_1(0)\cap \{u_{i,k}>0\} \rightarrow
B_1(0)\cap \{y\cdot\nu_i<0\}, \mbox{ as } \rho\to 0,\ \mbox{and}\\
& \nabla u_{i,k} \rightarrow \nabla u_i=-\lambda_{i} \nu_i
\chi_{\{y\cdot\nu_i<0\}} \mbox{ a.e in }B_1(0).
\end{align*}

Therefore
$$
\rho_k^{-N-1}\int_{B_{\rho_k}(x_i)}(|\nabla v_k|^p-|\nabla
u_{\ep}|^p)\, dx \rightarrow \int_{B_{1}(0)\cap \{y\cdot\nu_i>0\}}
\lambda_{i}^p \big(\mbox{div}(\eta_i)- p \,\nu_i^t\, D\eta_i
\,\nu_i\big)\, dy
$$

Using that
$$
\mbox{div}(\eta_i)- p\, \nu_i^t \,D\eta_i\, \nu_i = (-1)^i(1-p)
\frac{\phi'(|y|)}{|y|}( y\cdot\nu_i)= (-1)^i(1-p)
\mbox{div}(\eta_i),
$$
we obtain
$$
\rho_k^{-N-1}\int_{B_{\rho_k}(x_i)}(|\nabla v_k|^p-|\nabla
u_{\ep}|^p)\, dx \rightarrow (-1)^i(1-p) \lambda_{i}^p
\int_{B_{1}(0)\cap \{y\cdot\nu_i=0\}} \phi(|y|)\, d\H(y)
$$

Hence
\begin{equation}\label{estimacion.L2}
\begin{aligned}
& \int_{\Omega} |\nabla v_k|^p \, dx- \int_{\Omega} |\nabla
u_{\ep}|^p \, dx= \\
& \hskip1cm =\rho_k^{N+1} (\lambda_1^p-\lambda_0^p)
\int_{B_1(0)\cap\{y_1=0\}} (p-1) \phi(|y|)\, d\H(y)  +
o(\rho_k^{N+1}).
\end{aligned}
\end{equation}
If we take  $k$ large enough we get
$$
\J_\ep(v_k)< \J_\ep(u),
$$
a contradiction.
\end{proof}

\medskip

\begin{proof} [\bf Proof of Theorem \ref{const}]
Now, the Theorem follows as in steps 2 and 3 of Theorem 5.1 in
\cite{L}, using Lemma 5.4 in \cite{DP}, Theorem \ref{psub0} and
properties (1)--(8) of Lemma \ref{propblowup}. We sketch the proof
here for the reader's convenience.

Let $x_0\in\Omega\cap\partial\{u>0\}$ and let
$$
\lambda=\lambda(x_0):=\limsup_{\stackrel{x\to x_0}{u(x)>0}} |\nabla
u(x)|.
$$
Then there exists a sequence $z_k\rightarrow x_0$ such
that
$$
u(z_k)>0,\quad \quad |\nabla u(z_k)|\rightarrow \lambda.
$$
Let $y_k$ be the nearest point to $z_k$ on $\Omega \cap
\partial\{u>0\}$ and let $d_k = |z_k-y_k|$. Consider the blow up
sequence with respect to $B_{d_k}(y_k)$ with limit $u_0$, such
that there exists
$$
\nu:=\lim_{k\to\infty} \frac{y_k-z_k}{d_k},
$$
and suppose that $\nu = e_N$. Using the results of Appendix
\ref{appB}, we can proceed as in \cite{DP} p.13 to prove that
$0<\lambda<\infty$ and
$$
u_0(x)=-\lambda x_N \mbox{ in } \{x_N\leq 0\}.
$$
Finally by Lemma \ref{propblowup}(8) we have that $0\in
\partial \{u_0>0\}$ and then, using Lemma \ref{propblowup}(6) we
see that $u_0$ satisfies the hypotheses of Theorem \ref{psub0}.
Therefore $u_0=0$ in $\{x_N>0\}$. Then $u_0 = \lambda \max(-x\cdot
\nu,0)$.

To complete the proof, we follow the lines in step 3 of Theorem
5.1 in \cite{L}. This is, we apply Lemma \ref{blow} to this blow
up sequence and to a blow up sequence centered at a regular point
of the free boundary.

A similar argument proves \eqref{limsup.dist}.
\end{proof}

Summing up, we have the following theorem,
\begin{teo}\label{weak}
Let $u\in \K$ be a solution to \eqref{Pep}. Then $u$
is a weak solution to the following free boundary problem
\begin{align*}
& \Delta_p u = 0 \ \, \qquad \mbox{ in } \{u>0\}\cap\Omega,\\
& \frac{\partial u}{\partial\nu} = \lambda_{u} \, \qquad\mbox{ on
} \partial\{u>0\}\cap\Omega,
\end{align*}
where $\lambda_{u}$ is the constant in Theorem \ref{const}. More
precisely, $\H-$a.e. point $x_0\in\partial\{u>0\}$ belongs to
$\partial_{red}\{u>0\}$ and
$$
u(x_0+x) = \lambda_{u} (x\cdot \nu(x_0))^- + o(|x|) \quad
\mbox{for}\quad x\to 0.
$$
\end{teo}

Finally, we get an estimate of the gradient of $u$ that will be
needed in order to get the regularity of the free boundary.

\begin{teo}\label{estimacion.del.gradiente}
Let $u\in \K$ be a solution to \eqref{Pep}. Given $D \subset
\subset \Omega$, there exist constants $C=C(N,\ep,D)$ $r_0 =
r_0(N,D) > 0$ and $\gamma=\gamma(N,\ep,D)>0$ such that, if $x_0\in
D \cap \partial\{u>0\}$ and $r<r_0$, then
$$
\sup_{B_r(x_0)} |\nabla u|\le \lambda_u(1+Cr^\gamma).
$$
\end{teo}

\begin{proof}
The proof follows the lines  of the proof of Theorem 7.1 in
\cite{DP}.
\end{proof}

\medskip

As a corollary we have the following regularity result for the
free boundary $\partial \{u>0\}$.

\begin{corol}\label{teo.regularity}
Let $u_{\ep}\in \K$ be a solution to \eqref{Pep}. Then
$\partial_{red}\{u_{\ep}>0\}$ is a $C^{1,\beta}$ surface locally
in $\Omega$ and the remainder of the free boundary has
$\H-$measure zero. Moreover, if $N=2$ then the whole free boundary
is a $C^{1,\beta}$ surface.
\end{corol}

\begin{proof}
See \cite{DP} Corollary 9.2.
\end{proof}


\section{Behavior of the minimizer for small $\ep$.}
\label{sect.ep} \setcounter{equation}{0} In this section, since we
want to analyze the dependence of the problem with respect to
$\ep$ we will again denote by $\uep$ a solution to problem
\eqref{Pep}.

To complete the analysis of the problem, we will now show that if
$\ep$ is small enough, then
$$
|\{\uep>0\}| = \alpha.
$$
To this end, we need to prove that the constant $\lambda_\ep :=
\lambda_{\uep}$ is bounded from above and below by positive
constants independent of $\ep$. We perform this task in a series
of lemmas.

\begin{lema}\label{AAC-lema5}
Let $\uep\in \K$ be a solution to \eqref{Pep}. Then, there exists
a constant $C>0$ independent of $\ep$ such that
$$
\lambda_\ep := \lambda_{\uep} \le C.
$$
\end{lema}

\begin{proof}
The proof is similar to the one in  \cite{AAC}, Theorem 3.

First we will prove that there exist $C,c>0$, independent of
$\ep$, such that
$$
c\leq |\{\uep>0\}|\leq C\ep+\alpha.
$$
In fact, as in Theorem \ref{existencia} we have that
$F_{\ep}(|\{\uep>0\}|) \leq C$ thus obtaining the bound from
above. On the other hand, taking $q<p$, using the Sobolev trace
Theorem, the H\"{o}lder inequality and the fact that
$\|\uep\|_{\wp} \leq C$ (see Theorem \ref{existencia}) we have
$$
\int_{\partial \Omega} \varphi_0^q\, d\H \leq C
|\{\uep>0\}|^{\frac{p-q}{p}} \|u\|_{\wp}^q \leq C
|\{\uep>0\}|^{\frac{p-q}{p}},
$$
and thus we obtain the bound from below.

Take $D\subset\subset \Omega$ smooth, such that
$\theta=|D|>\alpha$ and $|\Omega\setminus D|<c$ then,
$$
|D\cap \{\uep>0\}|\leq\alpha +C\ep<\theta
$$
for $\ep$ small enough. On the other hand
$$
|D\cap  \{\uep>0\}|\geq |\{\uep>0\}|-|\Omega\setminus
D|\geq c-|\Omega\setminus D|>0,
$$
Therefore by the relative isoperimetric inequality we have
$$
\H(D\cap \partial\{\uep>0\})\geq c\min{\Big\{|D\cap\{\uep>0\}|,
|D\cap\{\uep=0\}|\Big\}}^{\frac{N-1}{N}}\geq c>0.
$$
Now let $w$ be the $p-$harmonic function in $\Omega$ with boundary
data equal to $\varphi_0$. Using Theorem \ref{DP-teo5.2} and
Theorem \ref{const} we have,
\begin{align*}
C & \geq\int_{\Omega} |\nabla \uep|^{p-2} \nabla \uep \nabla
(\uep-w)\, dx=\int_{\Omega\cap\partial\{\uep>0\}} w\lambda_{\ep}\,
d\H\geq \int_{D\cap\partial\{\uep>0\}}w\lambda_{\ep}\, d\H \\
&\geq \lambda_{\ep} (\inf_{D} w) \H(D\cap\partial\{\uep>0\})\geq c
\lambda_{\ep}.
\end{align*}
Now the result follows.
\end{proof}

\begin{lema}\label{gamapromedio}
Let $\uep\in \K$ be a solution to \eqref{Pep}, $B_r \subset\subset
\Omega$  and $v$ a solution to
$$
\Delta_p v=0\quad\mbox{in }B_r,\quad \quad v=u_{\ep}\quad\mbox{on
}\partial B_r.
$$
then
$$
\int_{B_r} |\nabla (u_{\ep}-v)|^q\, dx \geq C |B_r \cap
\{u_{\ep}=0\}| \left(\frac{1}{r}\left(\di\pint_{B_r}
u_{\ep}^{\gamma}\, dx \right)^{1/\gamma}\right)^q
$$
for all $q\geq 1$ and for any $\gamma<\frac{N(p-1)}{N-p}$ if
$p\leq N$, $\gamma<\infty$ if $p>N$, and $C$ is a constant
independent of $\ep$.
\end{lema}

\begin{proof}
The idea of the proof is similar to Lemma 3.2 in \cite{AC}. We
include  the details since there are differences due to the fact
that we are dealing with the $p$-laplacian instead of the
laplacian.

First let us assume that $B_r = B_1(0)$. For $|z|\le \frac12$ we
consider the change of variables from $B_1$ into itself such that
$z$ becomes the new origin. We call $u_z(x) = u \big( (1-|x|)z + x
\big)$, $v_z(x) = v\big((1-|x|)z+x\big)$ and define
$$
r_\xi = \inf\Big\{ r \,/\, \frac18\le r\le 1\quad\mbox{and}\quad u_z(r\xi)=0\Big\},
$$
if this set is nonempty. Observe that this change of variables
leaves the boundary fixed.

Now, for almost every $\xi\in \partial B_1$ we have
\begin{equation}\label{cota-sup}
v_z(r_\xi \xi) = \int_{r_\xi}^1 \frac{d}{dr}(u_z-v_z)(r\xi)\,
dr\le (1-r_\xi)^{1/q'}\left(\int_{r_\xi}^1
|\nabla(u_z-v_z)(r\xi)|^q\, dr\right)^{1/q}.
\end{equation}

Let us assume that the following inequality holds
\begin{equation}\label{v_z}
v_z(r_\xi \xi)\ge C(N,\Omega)(1-r_\xi)\left(\pint_{ B_1}
u^{\gamma}\, dx\right)^{1/\gamma}.
\end{equation}

Then, using \eqref{cota-sup} and \eqref{v_z}, integrating first
over $\partial B_1$ and then over $|z|\leq 1/2$ we obtain as in
\cite{AC},
$$
\int_{B_1} |\nabla (u-v)|^q\, dx \geq C |B_1 \cap \{u=0\}|
\left(\di\pint_{B_1} u^{\gamma}\, dx\right)^{q/\gamma}.
$$
If we take $u_r(x)=\frac{1}{r}u(x_0+rx)$ (where $x_0$ is the
center of the ball $B_r$) then
\begin{align*}
& \int_{B_1} |\nabla (u_r-v_r)|^q\, dx =r^{-N}
\int_{B_r} |\nabla (u-v)|^q\, dy,\\
& |B_1\cap \{u_r=0\}|=r^{-N}
|B_r\cap \{u=0\}| \quad\mbox{ and }\\
& \left(\di\pint_{B_1}u_r^{\gamma}\, dx\right)^{1/\gamma}=
\frac{1}{r}\left(\di\pint_{B_r} u^{\gamma}\, dy\right)^{1/\gamma},
\end{align*}
so we have the desired result.

Therefore we only have to prove \eqref{v_z}. If $|(1-r_\xi)z+r_\xi
\xi|\le \frac34$, by Harnack inequality,
$$
v_z(r_\xi \xi)\ge C_N v(0).
$$
By Theorem 1.2 in \cite{T} we have
\begin{equation}\label{truding}
v(0)\ge \alpha(N,\Omega) \left(\pint_{B_1} v^{\gamma}\,
dx\right)^{1/\gamma} \geq \alpha(N,\Omega) \left(\pint_{B_1}
u^{\gamma}\, dx\right)^{1/\gamma}.
\end{equation}

If $|(1-r_\xi)z+r_\xi \xi|\ge \frac34$ we prove by a comparison
argument that inequality \eqref{v_z} also holds. In fact, again by
Theorem 1.2 in \cite{T},
$$
v\ge C_N \alpha\displaystyle \left(\pint_{B_1}
u^{\gamma}\, dx\right)^{1/\gamma} \mbox{ in }B_{3/4}.
$$
Let $w(x) = e^{-\lambda|x|^2} - e^{-\lambda}$. There exists
$\lambda=\lambda(N,\alpha)$ such that
$$
\begin{cases}
\Delta_p w \ge 0 & \mbox{in } B_1\setminus B_{3/4},\\
w \le C_N \alpha & \mbox{on } \partial B_{3/4},\\
w=0 & \mbox{on } \partial B_1,
\end{cases}
$$
so that,
$$
v\ge w \displaystyle\left(\pint_{B_1} u^{\gamma}\,
dx\right)^{1/\gamma}\ge C (1-|x|) \left(\pint_{B_1} u^{\gamma}\,
dx\right)^{1/\gamma}\quad \mbox{in}\quad B_1\setminus B_{3/4}.
$$
Therefore
$$
v_z(r_\xi \xi)\ge C \Big(1-|(1-r_\xi)z+r_\xi \xi|\Big)
\left(\pint_{B_1} u^{\gamma}\, dx\right)^{1/\gamma}\ge C (1-r_\xi)
\left(\pint_{B_1} u^{\gamma}\, dx\right)^{1/\gamma}
$$
since $|z|\le \frac12$. So that \eqref{v_z} holds for every
$r_\xi\ge \frac18$.

This completes the proof.
\end{proof}

\medskip

\begin{lema}\label{cotaabajo}
Let $\uep\in\mathcal{K}$ be a solution to $(P_{\ep})$,
then
$$
\lambda_{\ep}\geq c>0,
$$
where $c$ is independent of $\ep$
\end{lema}

\begin{proof}
We proceed as in Lemma 6 in \cite{AAC}. We will use the following fact that we
prove in Lemma \ref{pmayor2} bellow:
For every $\ep>0$ there is a neighborhood of $A$ in $\Omega$ where $\uep>0$.

Let $y_0\in  A$ and let $D_t$ with $0\le t\le 1$ be a family
of open sets with smooth boundary and uniformly (in $\ep$ and $t$)
bounded curvatures such that $D_0$ is an exterior tangent ball at
$y_0$, $D_1$ contains a free boundary point, $D_0\subset\subset
D_t$ for $t>0$ and $D_t\cap\partial\Omega\subset  A$.

Let $t\in (0,1]$ be the first time such that $D_t$ touches the
free boundary and let $x_0\in \partial D_t\cap \partial \{\uep>0\}
\cap \Omega$. Now, take $w$ such that $\Delta_p w=0$ in $D_t
\setminus \overline{D}_0$ with $w=c_0$ on $\partial D_0$ and $w=0$
on $\partial D_t$. Thus $w \le \uep$ in $D_t\cap\Omega$ and
$\partial_{-\nu} w(x_0) \ge c\, c_0$ with $c>0$ independent of
$\ep$. Therefore, for $r$ small enough,
\begin{equation}\label{gamma2}
\left(\pint_{B_r(x_0)} \uep^{\gamma}\, dx\right)^{1/\gamma} \ge
\left(\pint_{ B_r(x_0)} w^{\gamma}\, dx\right)^{1/\gamma} \ge r
\bar c\, c_0,
\end{equation}
with $\bar c$ is independent of $\ep$.

If $v_0$ is the solution to
$$
\begin{cases}
\Delta_p v_0 = 0 & \mbox{in } B_r(x_0)\\
v_0 = \uep & \mbox{on } \partial B_r(x_0),
\end{cases}
$$
then, by Lemma \ref{gamapromedio}, we have
\begin{align*}
& \int_{B_r} |\nabla (u_{\ep}-v_0)|^p\, dx \geq C |B_r \cap
\{u_{\ep}=0\}| \left(\frac{1}{r}\left(\di\pint_{B_r}
u_{\ep}^{\gamma}\, dx \right)^{1/\gamma}\right)^p, \mbox{ for }
p\geq 2\\
& \int_{B_r} |\nabla (u_{\ep}-v_0)|^2 \, dx\geq C |B_r \cap
\{u_{\ep}=0\}| \left(\frac{1}{r}\left(\di\pint_{B_r}
u_{\ep}^{\gamma}\, dx\right)^{1/\gamma}\right)^2, \mbox{ for }
1<p\leq 2.
\end{align*}

Then using \eqref{desigualdadDP1} we obtain,
\begin{equation}\label{primera}
\int_{B_r} (|\nabla u_{\ep}|^p-|\nabla v_0|^p)\, dx \geq C |B_r
\cap \{u_{\ep}=0\}| \left(\frac{1}{r}\left(\di\pint_{B_r}
u_{\ep}^{\gamma}\, dx \right)^{1/\gamma}\right)^p
\end{equation}
for $p\ge2$.

By Theorem \ref{const} and Lemma \ref{AAC-lema5} we have that,
near $x_0$,  $|\nabla u_{\ep}|$ is bounded from above by a
constant independent of $\ep$. Then by \eqref{desigualdadDP2} we
obtain that \eqref{primera} also holds  for $1<p\le2$ if $r$ is
small enough (depending on $\ep$). Then by \eqref{gamma2}
\begin{equation}\label{cota1}
\int_{B_r(x_0)} (|\nabla u_{\ep}|^p-|\nabla v_0|^p)\, dx \geq c
\delta_r
\end{equation}
where $\delta_r= |B_r(x_0)\cap \{u_{\ep}=0\}|$ and $c$ is a
constant independent of $\ep$.

Consider now a free boundary point $x_1$ away from $x_0$. We can
choose $x_1$ to be regular.

Let us take
$$
\tau_{\rho}(x)=\begin{cases} \displaystyle
x-\rho^2\phi\left(\frac{|x-x_1|}{\rho}\right)
\nu_{u_{\ep}}(x_1) & \mbox{ for } x\in B_{\rho}(x_1),\\
x &\mbox{ elsewhere, }
\end{cases}
$$
where $\phi\in C_0^{\infty}(-1,1)$ with $\phi'(0)=0$.

Now choose $\rho$ such that
$$
\delta_r=\rho^2 \int_{B_{\rho}(x_1)\cap\partial\{u_{\ep}>0\}}
\phi\left(\frac{|x-x_1|}{\rho}\right)\, d\H.
$$
Take $v_\rho(\tau_{\rho}(x))=u_{\ep}(x)$ and
$$
v= \begin{cases}
v_0 & \mbox{ in } B_{r}(x_0)\\
v_{\rho} & \mbox{ in } B_{\rho}(x_1)\\
u_{\ep} & \mbox{ elsewhere.}
\end{cases}
$$
Thus, we have that
\begin{equation}\label{medida}
|\{v>0\}|=|\{u_{\ep}>0\}|.
\end{equation}
On the other hand as in Lemma \ref{blow}, we have
\begin{align*}
\int_{B_{\rho}(x_1)}(|\nabla v_\rho|^p-|\nabla u_{\ep}|^p)\ dy & =
\int_{\tau_{\rho}\left(B_{\rho}(x_1)\right)\cap
\{v_\rho>0\}}|\nabla v_\rho|^p\ dy-\int_{ B_{\rho}(x_1)}|\nabla
u|^p\ dx \\
& = \int_{B_{\rho}(x_1)\cap \{u>0\}}\rho (|\nabla u_{\ep}|^p
\mbox{div}\,\eta- p \,|\nabla u_{\ep}|^{p-2} \nabla u_{\ep} D\eta
\nabla u_{\ep}) + o(\rho)\ dx
\end{align*}
where $\eta(y)=-\phi(|y|) \nu(x_1)$.  Using the fact that $\eta$
is bounded from above by a constant $k$ independent of $\rho$ and
$\ep$, and that $|\nabla u_{\ep}| = \lambda_{\ep} + O(\rho^2)$ in
$B_{\rho}(x_1)$ we have
$$
\int_{B_{\rho}(x_1)}(|\nabla v_\rho|^p-|\nabla u_{\ep}|^p)\, dy
\leq k\lambda_{\ep}^p \rho^{N+1}+ o(\rho)\rho^N
$$
but, $\delta_r$ has the same order of $\rho^{N+1}$ then
\begin{equation}\label{cota2}
\int_{B_{\rho}(x_1)}(|\nabla v_\rho|^p-|\nabla u_{\ep}|^p)\, dy
\leq k\lambda_{\ep}^p \delta_r+ o(\delta_r).
\end{equation}
Therefore by \eqref{cota1},  \eqref{cota2} and \eqref{medida} we
have
$$
0\leq \jep(v)-\jep(u_{\ep})\leq -c\delta_r+k\lambda_{\ep}^p \delta_r+o(\delta_r)
$$
and then $\lambda_{\ep}\geq c>0$.
\end{proof}

Now we prove the positivity result that was used in the previous Lemma.
\begin{lema}\label{pmayor2}
For every $\ep>0$ there exists a neighborhood of $A$ in $\Omega$ such that
$\uep>0$ in this neighborhood.
\end{lema}
\begin{proof} Let $y_0\in A$ and let $B_\delta(z_0)$ be an exterior tangent ball to $\partial\Omega$
at $y_0$ such
that $\overline\Omega\cap \overline B=\{y_0\}$. Let us take $\delta$ small enough
so that $B_{2\delta}(z_0)\cap\partial \Omega\subset\subset A$. Let $\wep$ be a minimizer of
\begin{equation}\label{aux}
\widetilde J_\ep(w):= \int_\calr |\nabla w|^p\,dx + \frac 1\ep |\{w>0\}\cap\calr|
\end{equation}
in $\{w\in W^{1,p}(\calr)\,,\, w=0\mbox{ on }\partial B_{2\delta}(z_0),\ w=c_0
\mbox{ on }\partial B_{\delta}(z_0)\}$. Here $\calr= B_{2\delta}(z_0)\setminus \overline B_{\delta}(z_0)$.

Every minimizer of \eqref{aux} is radially symmetric and radially decreasing with
respect to $z_0$. This is seen by using  Schwartz symmetrization after extending $\wep$ to $B_\delta(z_0)$
as the constant function $c_0$ (see \cite{K}). This symmetrization
preserves the distribution function and strictly decreases the $L^p$ norm of the gradient unless
the function is already radially symmetric and radially decreasing. Moreover, these minimizers are
ordered and their supports are nested. Let us take as $\wep$ the smallest minimizer.

By the properties of $\wep$ there holds that $\wep$ is strictly positive in a ring around $B_{\delta}(z_0)$.
Also $\wep$ is continuous in $\calr$. Recall that $\uep$ is continuous in $\Omega$.
Let us see that $\uep\ge\wep$ in $\calr\cap\Omega$. This will prove the statement.

Assume instead that $\{\wep>\uep\}\neq\emptyset$.

Let us first consider the function $ v=\min\{\uep,\wep\}$ in $\calr\cap\Omega$. Since $\uep\ge c_0\ge \wep$ on
$\partial\Omega\cap\calr$ and $\uep\ge0=\wep$ on $\Omega\cap\partial\calr$ there holds that $v=\wep$ on
$\partial(\calr\cap\Omega)$. Therefore, the function $\underline v=v$ in $\calr\cap\Omega$, $\underline v=\wep$
in $\calr\setminus\Omega$ is an admissible function for the minimization problem \eqref{aux}. Since $\wep$ is
the smallest minimizer and, by assumption $\underline v\le \wep$ and $\underline v\neq \wep$, there holds that $\widetilde J_\ep(\underline v)>
\widetilde J_\ep(\wep)$. Since $\underline v=\wep$ in $\calr\setminus\Omega$ and in $\calr\cap\Omega\cap\{\wep\le\uep\}$
and equal to $\uep$ outside those sets
there holds that (with $\cald=\calr\cap\Omega\cap\{\wep>\uep\}$),
\begin{equation}\label{aux1}
\int_{\cald} |\nabla \uep|^p\,dx +\frac 1\ep |\{\uep>0\}\cap \cald|>
\int_{\cald} |\nabla \wep|^p\,dx +\frac 1\ep |\{\wep>0\}\cap \cald|.
\end{equation}

Let now $\bar v=\max\{\uep,\wep\}$ in $\calr\cap\Omega$, $\bar v=\uep$ in $\Omega\setminus \calr$. This function
is admissible for \eqref{Pep} so that
$$
\int_\Omega|\nabla\bar v|^p\,dx +F_\ep\big(|\{\bar v>0\}|\big)\ge
\int_\Omega|\nabla\uep|^p\,dx +F_\ep\big(|\{\uep>0\}|\big).
$$

Since $\bar v=\wep$ in $\cald$ and $\bar v=\uep$ in $\Omega\setminus\cald$,
\begin{equation}\label{aux2}
\begin{aligned}
&\int_\cald|\nabla\wep|^p\,dx +F_\ep\big(|\{\uep>0\}|+|\{\wep>0\}\cap\cald|-|\{\uep>0\}\cap\cald|\big)\\
&\hskip1.5cm\ge\int_\cald|\nabla\uep|^p\,dx +F_\ep\big(|\{\uep>0\}|\big).
\end{aligned}
\end{equation}

By \eqref{aux1} and \eqref{aux2} (with $C_w=|\{\wep>0\}\cap\cald|$ and $C_u=|\{\uep>0\}\cap\cald|$) we have,
$$
\begin{aligned}
&\int_\cald|\nabla\uep|^p\,dx>\int_\cald|\nabla \wep|^p\,dx +\frac1\ep (C_w-C_u)\\
\ge
&\int_\cald|\nabla\uep|^p\,dx+F_\ep(|\{\uep>0\}|)-F_\ep(|\{\uep>0\}|+C_w-C_u)+\frac1\ep(C_w-C_u).
\end{aligned}
$$
Thus,
$$
F_\ep(|\{\uep>0\}|+C_w-C_u)-F_\ep(|\{\uep>0\})>\frac1\ep(C_w-C_u)
$$
which is a contradiction since $F_\ep(A)-F_\ep(B)\le \frac1\ep(A-B)$ if $A\ge B$ and $C_w\ge C_u$ by assumption.

Therefore, $\uep\ge\wep$ in $\calr\cap\Omega$ and the lemma is proved.
\end{proof}

With these uniform bounds on $\lambda_{\ep}$, we can prove the desired result.
\begin{teo}\label{final}
There exists
$\ep_0>0$ such that if $\uep\in \K$ is a solution to \eqref{Pep} and  $\ep<\ep_0$ there holds that
 $|\{u_{\ep}>0\}|=\alpha$.
Therefore, $u_{\ep}$ is a minimizer of $\J$ in $\K_{\alpha}$.
\end{teo}
\begin{proof}
Arguing by contradiction, assume first that $|\{\uep > 0\}| >
\alpha$. Let $x_1 \in \partial \{ \uep >0 \} \cap \Omega$ be a
regular point. We will proceed as in the proof of Lemma
\ref{cotaabajo}. Given $\delta >0$, we perturb the domain $\{\uep
>0\}$ in a neighborhood of $x_1$, decreasing its measure by
$\delta$. We choose $\delta$ small so that the measure of the
perturbed set is still larger than $\alpha$. Take
$v_{\rho}(\tau_{\rho}(x)) = u_{\ep}(x)$, and  let
$$
v=\begin{cases}
v_{\rho} & \mbox{ in } B_{\rho}(x_1)\\
u_{\ep} & \mbox{ elsewhere,}
\end{cases}
$$
where $\tau_{\rho}$ is the function that we have considered in the
previous lemma.

We have
\begin{align*}
0 \le \J_\ep(v) - \J_\ep(\uep) & = \int_\Omega |\nabla v|^p\, dx -
\int_\Omega |\nabla \uep|^p\, dx+ F_\ep (| \{ v>0\}|)
- F_\ep (| \{ \uep >0\}|) \\
& \le k\lambda_\ep^p\delta + o_\ep(\delta) - \frac{1}{\ep} \delta
\le \left(kC^p - \frac{1}{\ep}\right) \delta +  o_\ep(\delta) < 0,
\end{align*}
if $\ep < \ep_0$ and then $\delta < \delta_0 (\ep)$. A
contradiction.

Now assume that $|\{\uep > 0\}| < \alpha$. This case, is a little
bit different from the other. First, we proceed as in the previous
case but this time we perturb in a neighborhood of $x_1$ the set
$\{\uep >0\}$ increasing its measure by $\delta$. That is, take
$$
\tau_{\rho}(x) = \begin{cases} \displaystyle
x+\rho^2\phi\left(\frac{|x-x_1|}{\rho}\right)
\nu_{u_{\ep}}(x_1) & \mbox{ for } x\in B_{\rho}(x_1),\\
x & \mbox{ elsewhere,}
\end{cases}
$$
where $\phi\in C_0^{\infty}$ supported in the unit interval, take
$v_\rho(\tau_{\rho}(x))=u_{\ep}(x)$ and
$$
v = \begin{cases}
v_{\rho} & \mbox{ in } B_{\rho}(x_1)\\
u_{\ep} & \mbox{ elsewhere.}
\end{cases}
$$
For $\rho$ small enough  we have $|\{v>0\}|<\alpha$
and
$$
|\{v>0\}|-|\{u_{\ep}>0\}|=C\rho^{N+1}+o(\rho^{N+1}),
$$
therefore
\begin{equation}\label{chica1}
F_{\ep}(|\{v>0\}|)-F_{\ep}(|\{u_{\ep}>0\}|)\leq C\ep
\rho^{N+1}+o_{\ep}(\rho^{N+1}).
\end{equation}

In order to estimate the other term,  we will make use of a blow
up argument as in Lemma \ref{blow}. In fact, we take $u_{\rho}(y)
= \frac{1}{\rho}u(x_1+\rho y)$ and we change variables to obtain,
\begin{align*}
\rho^{-N}\int_{B_{\rho}(x_1)}&(|\nabla v_\rho|^p-|\nabla
u_{\ep}|^p)\, dx \\
&= \int_{B_{1}(0)\cap \{u_{\rho}>0\}}\rho [|\nabla u_{\rho}|^p
\mbox{div}(\eta)- p |\nabla u_{\rho}|^{p-2} (\nabla u_{\rho})^t
D\eta \nabla u_{\rho}] + o(\rho)\, dy
\end{align*}
where $\eta(y)=\phi(|y|) \nu(x_1)$. Now, as in Lemma \ref{blow} we
get,
$$
\rho^{-N-1}\int_{B_{\rho}(x_1)}(|\nabla v_\rho|^p-|\nabla
u_{\ep}|^p)\, dx \rightarrow (1-p)
\lambda_{\ep}^p\int_{B_{1}(0)\cap \{y\cdot\nu=0\}}\phi(|y|)\,
d\H(y).
$$
Therefore
\begin{align}\label{chica2}
\int_{B_{\rho}(x_1)}(|\nabla v_\rho|^p - |\nabla u_{\ep}|^p)\, dx
= C \rho^{N+1}(1-p) \lambda_{\ep}^p + o(\rho^{N+1}).
\end{align}
Finally, combining \eqref{chica1} and \eqref{chica2} we have
\begin{align*}
0 \le \J_\ep(v ) - \J_\ep(\uep) & = \int_\Omega |\nabla v|^p\, dx
- \int_\Omega |\nabla \uep|^p\, dx + F_\ep (| \{ v>0\}|)
- F_\ep (| \{ \uep >0\}|) \\
& \le C(1-p)\lambda_\ep^p\delta + o_\ep(\delta) + C\ep \delta \le
C(-c^p + \ep) \delta +  o_\ep(\delta) < 0,
\end{align*}
if $\ep < \ep_1$ and then $\delta < \delta_0 (\ep)$. Again a
contradiction that ends the proof.
\end{proof}

As a corollary, we have the desired result for our problem

\begin{corol}
For $\ep$ small
any minimizer  $u$ of $\J$ in $\K_{\alpha}$ is
a locally Lipschitz continuous function and
$\partial_{red}\{u>0\}$ is a $C^{1,\beta}$ surface locally in
$\Omega$ and the remainder of the free boundary has $\H-$measure
zero. In particular, this is the case for $p\ge 2$.
\end{corol}

\begin{proof}
If $u$ is minimizer of $\J$ in $\K_{\alpha}$, by  Theorem
\ref{final} we have that for small $\ep$ there exists a solution
$\uep$ to $(P_{\ep})$ such that $|\{\uep>0\}|=\alpha$, then $u$ is
a solution to $(P_{\ep})$, therefore the result follows.
\end{proof}

\medskip
\appendix
\renewcommand{\theequation}{\Alph{section}.\arabic{equation}}
\section{A result on $p$-harmonic functions with linear growth}\label{appA}
\setcounter{equation}{0}

In this section we will prove some properties of $p$-subharmonic
functions. From now on, we note $B_r^+= B_{r}(0)\cap\{x_N>0\}$.

\begin{teo}\label{psub0}
Let $u$ be a Lipschitz function in $\RR^N$ such that
\begin{enumerate}
\item $u\geq 0$ in $\RR^N$, $\Delta_p u=0$ in $\{u>0\}$.

\item $\{x_N<0\}\subset \{u>0\}$, $u=0$ in $\{x_N=0\}$.

\item There exists $0<\lambda_0<1$ such that
$\displaystyle\frac{|\{u=0\}\cap B_R(0)|}{|B_R(0)|}>\lambda_0$,
$\forall R>0$.
\end{enumerate}
Then $u=0$ in $\{x_N>0\}$.
\end{teo}

In order to prove this theorem we follow ideas from \cite{L}. To
this end, we need to prove a couple of lemmas.

\begin{lema}\label{psub2}
Let $u$ be a $p-$subharmonic function in $B_r^+$ such that, $0\leq
u\leq \alpha x_N$ in $B_r^+$, $u\leq \delta_0 \alpha x_N$ on
$\partial B_r^+\cap B_{r_0}(\bar{x})$ with $\bar{x}\in
\partial B^+_r$, $\bar{x}_N>0$ and $0<\delta_0<1$.

Then there exists $0<\gamma<1$ and $0<\ep\leq 1$, depending only
on $r$ and $N$, such that $u(x)\leq \gamma \alpha x_N$ in
$B_{\ep}^+$.
\end{lema}

\begin{proof}
By homogeneity of the $p-$laplacian we can suppose that $r=1$.

Let $\psi$ be a $p$-harmonic function in $B_1^+$, with smooth
boundary data, such that
$$
\begin{cases}
\psi=x_N & \mbox{on } \partial B_1^+\setminus B_{r_0}(\bar{x})\\
\delta_0 x_N\leq \psi \leq x_N & \mbox{on } \partial
B_1^+\cap B_{r_0}(\bar{x})\\
\psi=\delta_0 x_N & \mbox{on } \partial B_1^+\cap
B_{r_0/2}(\bar{x}).
\end{cases}
$$
Therefore, by comparison $u\leq \alpha \psi$ in $B_1^+$. Let us
see that there exist $0<\gamma<1$ and $\ep>0$, independent of
$\alpha$, such that $\psi\leq \gamma x_N$ in $B_{\ep}^+$.

First, $\psi\in C^{1,\beta}(\overline{B_1^+})$ for some $\beta>0$.
Then, (cf. \cite{JLM}) $\psi$ is a viscosity solution of
$$
|\nabla \psi|^{p-4}\Big\{|\nabla \psi|^2 \Delta \psi +(p-2)
\sum_{i,j=1}^N \psi_{x_i} \psi_{x_j} \psi_{x_i x_j} \Big\}=0.
$$
If $|\nabla \psi|\geq \mu>0$ in some open set $U$, we have that
$\psi$ is a solution of the linear uniformly elliptic equation
\begin{equation}\label{lue}
\sum_{i,j=1}^N a_{ij} \psi_{x_i x_j} =0 \quad \mbox{ in } U,
\end{equation}
where
$$
\min\{1,p-1\} |\nabla \psi|^{2} |\xi|^2\leq \sum_{i,j}^N
a_{ij} \xi_i \xi_j \leq \max\{1,p-1\} |\nabla \psi|^{2} |\xi|^2.
$$
Therefore, $\psi\in C^{2,\beta}(U)$ and is a classic solution of
\eqref{lue}.

Let $w=x_N-\psi$ then $w\in C^{1,\beta}(\overline{B_1^+})$ and is
a solution of
$$\mathcal{L}w=\sum_{i,j=1}^N a_{ij} w_{x_i x_j} =0$$
in any open set $U$ where $|\nabla \psi|\geq \mu>0$.

On the other hand, as $\psi\leq
x_N$ in $\partial B_1^+$ and both functions are solutions of the
$p-$laplacian we have, by comparison, that $\psi\leq x_N$ in
$B_1^+$. Therefore $w\geq 0$ in $B_1^+$.

Moreover, we have $w>0$ in $B_1^+$. In fact, suppose that there
exists $x_0$ such that $\psi(x_0)=x_{0,N}$. As $\psi\leq x_N$, we
have that $\nabla \psi (x_0)=e_N$. Then $|\nabla \psi (x_0)|=1$
and by continuity, $|\nabla \psi|\geq \frac{1}{2} >0$ in a
neighborhood $U$ of $x_0$. Therefore $w\geq 0$, $w(x_0)=0$ and
$\mathcal{L}w=0$ in $U$ with $\mathcal{L}$ uniformly elliptic in
$U$. Then by the strong maximum principle, $w\equiv 0$ in $U$.

So, we have that the set
$$
\mathcal{A}=\{x\in B_1^+\ /\ \psi(x)=x_N\},
$$
is a relative open and close subset of $B_1^+$. Then if there
exists $x_0$ such that $\psi(x_0)=x_{0,N}$, we have that
$\psi\equiv x_N$. Since this is not the case in some part of
$\partial B_1^+$, we arrive at a contradiction. Therefore
$\psi<x_N$, and this implies that $w>0$.

On the other hand, since $\psi\le x_N$ and $\psi=0$ on
$B_1^+\cap\{x_N=0\}$, we have that $\psi_{x_N} \leq 1$ on
$B_1^+\cap\{x_N=0\}$. Let us see that $\psi_{x_N} < 1$ in
$B_1^+\cap \{x_N=0\}$.

Assume that there exists $x_0\in B_1^+\cap \{x_N=0\}$ such that
$\psi_{x_N}(x_0)=1$ (so that $w_{x_N}(x_0)=0$). Then, $|\nabla
\psi|\ge 1/2$ in a neighborhood of $x_0$.  But $w$ is a positive
solution of $\mathcal{L}w=0$ in $B_1^+\cap B_{r_0}(x_0)$ for some
$r_0>0$, with $\mathcal{L}$ uniformly elliptic and $w=0$ on
$\{x_N=0\}$. Thus, by Hopf's Lemma, $w_{x_N}(x_0)>0$, a
contradiction.

Therefore $\psi_{x_N}<1$ in $B_1^+\cap\{x_N=0\}$. This implies
that there exists $0<\gamma<1$ and $\ep>0$ such that $\psi_{x_N} <
\gamma$ in $B_{\ep}^+$. From this, $\psi\leq \gamma x_N$ in
$B_{\ep}^+$, and then we have $u\leq \gamma \alpha x_N$ in
$B_{\ep}^+$, where $\ep$ and $\gamma$ only depend on $\psi$.
\end{proof}

\begin{lema}\label{psub}
Let $w$ be a function that satisfies,
\begin{enumerate}
\item $w$ is a Lipschitz function in $\RR^N$ with constant $L$.

\item $w\geq 0$ in $\RR^N$, $\Delta_p w=0$ in $\{w>0\}$.

\item $\{x_N<0\}\subset \{w>0\}$, $w=0$ in $\{x_N=0\}$.

\item There exists $0<\lambda_0<1$ such that
$\displaystyle\frac{|\{w=0\}\cap B_1(0)|}{|B_1(0)|}>\lambda_0$.

\item There exists $0\leq \alpha\leq L$ such that $w(x)\leq \alpha
x_N$ in $B_1(0)\cap\{x_N>0\}$.
\end{enumerate}
Then there exists $0<\gamma<1$ and $0<\ep\leq 1$ depending only on
$\lambda_0$ and $N$, such that $w(x)\leq \gamma \alpha x_N$ in
$B_{\ep}(0)\cap\{x_N>0\}$.
\end{lema}

\begin{proof}
Let $\beta=\frac{\lambda_0}{2^{N-1}}<1$, then by (3) and (4) there
exists $x_0\in B_1(0)$, with $x_{0,N}>\beta$ such that $w(x_0)=0$.
By (1), $w(x)\leq L|x-x_0|$, then if we take $r_0=\frac{\alpha
\beta}{4L}$, we have $w(x)\leq \frac{\alpha\beta}{4}$ for
$|x-x_0|<r_0$. As $\alpha/L\leq 1$, in that set there holds that
$x_N\geq \frac{3\beta}{4}$. Then we have that
$$
w(x)\leq \frac{\alpha x_N}{3} \mbox{ in } \partial B^+_r\cap
B_{r_0}(x_0),
$$
where $r=|x_0|>\beta$. Taking in Lemma \ref{psub2} $\delta_0=1/3$
and $\bar{x}=x_0$ we have that there exists $0<\gamma<1$ and
$0<\ep\leq 1$, depending on $r$ and $N$, such that $w(x)\leq
\gamma \alpha x_N$ in $B_{\ep}^+$.

As $r>\beta$ what we obtain is that $\gamma$ and $\ep$ only depend
on $\lambda_0$. Therefore the result follows.
\end{proof}

Now we are ready to proceed with the proof of the theorem.

\begin{proof} [\bf Proof of Theorem \ref{psub0}]
Once we have proved Lemma \ref{psub} we consider the same
iteration as in Theorem A.1, step 2 in \cite{L} and the result
follows.
\end{proof}

As a remark we mention that with Lemma \ref{psub2} we can also
prove the asymptotic development of $p-$harmonic functions, that
is
\begin{lema}
Let $u$ be  Lipschitz continuous in $\overline{B_1^+}$, $u\geq 0$
in $B_1^+$, $p-$harmonic in $\{u>0\}$ and vanishing on $\partial
B_1^+\cap \{x_N=0\}$. Then, in $B_1^+$, $u$ has the asymptotic
development
$$
u(x)=\alpha x_N+ o(|x|),
$$
with $\alpha\ge 0$.
\end{lema}

\begin{proof}
Let
$$
\alpha_j=\inf\{l\ /\ u\leq l x_n \mbox{ in } B_{2^{-j}}^+\}.
$$
Let $\alpha=\lim_{j\to\infty} \alpha_j$.

Given $\ep_0>0$ there exists $j_0$ such that for $j\geq j_0$ we
have $\alpha_j\leq \alpha+\ep_0$. From here, we have $u(x)\leq
(\alpha+\ep_0)x_N$ in $B_{2^{-k}}^+$ so that
$$
u(x)\leq \alpha x_N+o(|x|) \mbox{ in }B_1^+.
$$

If $\alpha=0$ the result follows. Assume that $\alpha>0$ and let
us suppose that $u(x)\neq \alpha x_N+o(|x|)$. Then there exists
$x_k\rightarrow 0$ and $\bar{\delta}>0$ such that
$$
u(x_k)\leq \alpha x_{k,N}-\bar{\delta} |x_k|.
$$
Let $r_k=|x_k|$ and $u_k(x)= r_k^{-1}u(r_k x)$. Then, there exists
$u_0$ such that, for a subsequence that we still call $u_k$,
$u_k\rightarrow u_0$ uniformly in $\overline{B_1^+}$ and
\begin{align*}
& u_k(\bar{x}_k)\leq \alpha \bar{x}_{k,N}-\bar{\delta}\\
& u_k(x)\leq (\alpha+\ep_0) x_N \mbox{ in } B_1^+,
\end{align*}
where $\bar{x}_k=\frac{x_k}{r_k}$, and we can assume that
$\bar x_k\rightarrow x_0$.

In fact, $u(x)\leq (\alpha+\ep_0) x_N$ in $B_{2^{-j_0}}^+$,
therefore $u_k(x)\leq (\alpha+\ep_0) x_N$ in $B_{r_k\,/\,
2^{-j_0}}^+$, and if $k$ is big enough $r_k\,/\, 2^{-j_0}\geq 1$.

If we take $\bar{\alpha}=\alpha+\ep_0$ we have
$$
\begin{cases}
\Delta_p u_k\geq 0 & \mbox{in }B_1^+\\
u_k=0 & \mbox{on } \{x_N=0\}\\
0\leq u_k \leq \bar{\alpha} x_N & \mbox{on } \partial
B_1^+\\
u_k\leq \delta_0 \bar{\alpha} x_N & \mbox{on } \partial B_1^+\cap
B_{\bar r}(\bar{x}),
\end{cases}
$$
for some $\bar{x} \in \partial B_1^+$, $\bar{x}_N>0$ and some
small $\bar r>0$.

In fact, as $u_k$ are continuous with uniform modulus of
continuity, we have
$$
u_k(x_0)\leq \alpha x_{0,N}-\frac{\bar{\delta}}{2}, \mbox{ if }
k\geq \bar{k}.
$$
Moreover there exists $r_0>0$ such that $u_k(x)\leq \alpha
x_N-\frac{\bar{\delta}}{4}$ in $B_{2r_0} (x_0)$. If $x_{0,N}>0$ we
take $\bar{x}=x_0$, if not, we take $\bar{x}\in B_{2r_0}(x_0)$
with $\bar{x}_N>0$ and
$$
u_k(x)\leq \alpha x_N-\frac{\bar{\delta}}{4}, \mbox{ in }
B_{r_0}(\bar{x})\subset\subset \{x_N>0\}.
$$
As $B_{r_0}(\bar{x})\subset\subset \{x_N>0\}$ there exists
$\delta_0$ such that $\alpha x_N-\frac{\bar{\delta}}{4}\leq
\delta_0 \alpha x_N \leq \delta_0 \bar{\alpha} x_N$ in
$B_{\bar r}(\bar{x})$ for some small $\bar r$, and the claim follows.

Now, by Lemma \ref{psub2}, there exists $0<\gamma<1$, $\ep>0$
independent of $\ep_0$ and $k$, such that $u_k(x)\leq \gamma
(\alpha+\ep_0) x_N$ in $B_{\ep}^+$. As $\gamma$ and $\ep$ are
independent of $k$ and $\ep_0$, taking $\ep_0\rightarrow 0$, we
have
$$
u_k(x)\leq \gamma \alpha x_N \mbox{ in }B_{\ep}^+.
$$
So that,
$$
u(x)\le \gamma\alpha x_N\mbox{ in }B_{r_k\ep}^+.
$$
Now if $j$ is big enough we have $\gamma \alpha <\alpha_j$ and
$2^{-j}\le r_k\ep$.  But this contradicts the definition of
$\alpha_j$. Therefore,
$$
u(x)=\alpha x_N+o(|x|),
$$
as we wanted to prove.
\end{proof}

\section{Blow-up limits}\label{appB}
\setcounter{equation}{0}

Now we  give the definition of blow-up sequence, and we collect
some properties of the limits of these blow-up sequences for
certain classes of functions that are used throughout the paper.

Let $u$ be a function with the following properties,
\begin{enumerate}
\item[(B1)] $u$ is Lipschitz in $\Omega$ with constant $L>0$,
$u\geq 0 \mbox{ in } \Omega$ and $\Delta_p u=0 \mbox{ in }
\Omega\cap\{u>0\}$.
\item[(B2)] Given $0<\kappa<1$, there exist
two positive constants $C_{\kappa}$ and $r_{\kappa}$ such that for
every ball $B_r(x_0)\subset\Omega$ and $0<r<r_{\kappa}$,
$$
\frac{1}{r}\left(\pint_{B_r(x_0)} u^{\gamma}\, dx
\right)^{1/\gamma}\leq C_{\kappa} \mbox{ implies that } u\equiv 0
\mbox{ in } B_{\kappa r}(x_0).
$$
\item[(B3)] There exist constants $r_0>0$ and
$0<\lambda_1\leq\lambda_2<1$ such that, for every ball
$B_r(x_0)\subset\Omega$ $x_0\mbox{ on }
\partial\{u>0\}$ and $0<r<r_0$
$$
\lambda_1\le
\frac{\left|B_r(x_0)\cap\{u>0\}\right|}{\left|B_r(x_0)\right|} \le
\lambda_2.
$$
\end{enumerate}

\begin{defi}
Let $B_{\rho_k}(x_k)\subset\Omega$ be a sequence of balls
with $\rho_k\to 0$, $x_k\to x_0\in \Omega$ and $u(x_k)=0$. Let
$$
u_k(x):=\frac{1}{\rho_k} u(x_k+\rho_k x).
$$
We call $u_k$ a blow-up sequence with respect to
$B_{\rho_k}(x_k)$.
\end{defi}

Since $u$ is locally Lipschitz continuous, there exists a blow-up
limit $u_0:\R^N\to\R$ such that for a subsequence,
\begin{align*}
& u_k\to u_0 \quad \mbox{in} \quad C^\alpha_{\rm loc}(\R^N)\quad
\mbox{for every}\quad 0<\alpha<1,\\
& \nabla u_k\to\nabla u_0\quad *-\mbox{weakly  in}\quad
L^\infty_{\rm loc}(\R^N),
\end{align*}
and $u_0$ is Lipschitz in $\RR^N$ with constant $L$.
\begin{lema}\label{propblowup}
If $u$ satisfies properties {\rm (B1), (B2)} and {\rm (B3)} then,
\begin{enumerate}
\item $u_0\geq 0$ in $\Omega$ and $\Delta_p u_0=0$ in $\{u_0>0\}$

\medskip

\item $\partial\{u_k>0\}\to \partial\{u_0>0\}$ locally in
Hausdorff distance,

\medskip

\item $\chi_{\{u_k>0\}}\to \chi_{\{u_0>0\}}$ in $L^1_{\rm
loc}(\R^N)$,

\medskip

\item If $K\subset\subset \{u_0=0\}$, then $u_k=0$ in $K$ for big
enough $k$,

\medskip

\item If $K\subset\subset \{u_0>0\}\cup \{u_0=0\}^\circ$, then
$\nabla u_k\rightarrow\nabla u_0$ uniformly in $K$,

\medskip

\item There exists a constant $0<\lambda<1$ such that,
$$
\frac{\left|B_R(y_0)\cap\{u_0=0\}\right|}{\left|B_R(y_0)\right|}\geq
\lambda, \quad \forall R>0, \forall y_0\in \partial\{u_0>0\}
$$

\medskip

\item $\nabla u_k\to\nabla u_0$ a.e in $\Omega$,

\medskip

\item If $x_k\in \partial\{u>0\}$, then $0\in
\partial\{u_0>0\}$
\end{enumerate}
\end{lema}

\begin{proof}
As $u_k$ are $p$-harmonic and $u_k\rightarrow u_0$ uniformly in
compacts subsets of $\RR^N$ then (1) holds. For the proof of
(2)--(8) see \cite{L}.
\end{proof}

\begin{ack} The authors want to thank Professor Arshak Petrosyan for providing the proof of Lemma \ref{pmayor2}.
\end{ack}

\end{document}